\documentclass[11pt,twoside,reqno]{amsart}
\usepackage[utf8]{inputenc}
\usepackage[T1]{fontenc}
\usepackage{import}
\usepackage{newclude}
\usepackage{anysize}
\usepackage{amsmath}
\usepackage{amsthm}
\usepackage{amsmath,amscd}
\usepackage{amssymb}
\usepackage{stmaryrd}
\usepackage{wasysym}
\usepackage{mathrsfs}
\usepackage{hyperref}
\usepackage{cleveref}
\usepackage{enumitem}
\usepackage{dsfont}
\usepackage{soul}
\hypersetup{colorlinks=black,linkcolor=black,citecolor=black}
\usepackage{color}
\usepackage[all]{xy}
\usepackage{float}
\usepackage{bm}
\usepackage{mathtools}
\usepackage{makecell}
\usepackage{thmtools}
\usepackage{setspace}
\usepackage{comment}
\usepackage{tikz}
\usepackage{tikz-cd}
\usepackage{mathdots}
\usepackage{graphicx}
\usepackage[myheadings]{fullpage}

\numberwithin{equation}{section}
\newcommand\mtop{.95in}
\newcommand\mbottom{.95in}
\newcommand\mleft{1in}
\newcommand\mright{1in}
\usepackage[top = \mtop, bottom = \mbottom, left = \mleft, right=\mright]{geometry}

\newtheorem{thm}{Theorem}[section]

\newtheorem{prop}[thm]{Proposition}

\newtheorem{lemma}[thm]{Lemma}

\theoremstyle{definition}

\newtheorem{rmk}[thm]{Remark}

\newcommand\reallywidehat[1]{%
\savestack{\tmpbox}{\stretchto{%
  \scaleto{%
    \scalerel*[\widthof{\ensuremath{#1}}]{\kern-.6pt\bigwedge\kern-.6pt}%
    {\rule[-\textheight/2]{1ex}{\textheight}}
  }{\textheight}%
}{0.5ex}}%
\stackon[1pt]{#1}{\tmpbox}%
}
\DeclareSymbolFont{bbold}{U}{bbold}{m}{n}
\DeclareSymbolFontAlphabet{\mathbbold}{bbold}

\makeatletter
\def\@tocline#1#2#3#4#5#6#7{\relax
  \ifnum #1>\c@tocdepth 
  \else
    \par \addpenalty\@secpenalty\addvspace{#2}%
    \begingroup \hyphenpenalty\@M
    \@ifempty{#4}{%
      \@tempdima\csname r@tocindent\number#1\endcsname\relax
    }{%
      \@tempdima#4\relax
    }%
    \parindent\z@ \leftskip#3\relax \advance\leftskip\@tempdima\relax
    \rightskip\@pnumwidth plus4em \parfillskip-\@pnumwidth
    #5\leavevmode\hskip-\@tempdima
      \ifcase #1
       \or\or \hskip 1em \or \hskip 2em \else \hskip 3em \fi%
      #6\nobreak\relax
    \hfill\hbox to\@pnumwidth{\@tocpagenum{#7}}\par
    \nobreak
    \endgroup
  \fi}
\makeatother

\newcommand{\R}{\mathbb{R}}

\newcommand{\Q}{\mathbb{Q}}

\newcommand{\C}{\mathbb{C}}
\newcommand{\F}{\mathbb{F}}

\newcommand{\I}{\mathcal{I}}
\newcommand{\mP}{\mathcal{P}}
\newcommand{\mL}{\mathcal{L}}

\newcommand{\mC}{\mathcal{C}}

\renewcommand{\char}{\text{char}}

\newcommand{\K}{K}

\title{The Szemerédi-Trotter theorem over arbitrary field of characteristic zero}
\author{Jiahe Shen}

\date{\today}
\begin{document}

\thanks{I thank my advisor Ivan Corwin for providing me funding support with his NSF grant DMS-2246576 and Simons Investigator grant 929852; Roger Van Peski, for carefully reading the draft and pointing out typos; Professor Ruixiang Zhang and Mehtaab Sawhney, for additional helpful conversations.}

\maketitle

\begin{abstract} 
Let $\mP$ be a set of $m$ points and $\mL$ a set of $n$ lines in $K^2$, where $K$ is a field with $\char(K)=0$. We prove the incidence bound
$$\I(\mP,\mL)=O(m^{2/3}n^{2/3}+m+n).$$
Moreover, this bound is sharp and cannot be improved. This resolves the Szemerédi-Trotter incidence problem for arbitrary field of characteristic zero.

The key tool of our proof is the Baby Lefschetz principle, which allows us to reduce the problem to the complex case. Based on this observation, we further derive several related results over $K$, including Beck's theorem, the Erd\H{o}s-Szemerédi sum–product estimate, and incidence theorems involving more general algebraic objects.
\end{abstract}



\tableofcontents

\section{Introduction}\label{sec: Introduction}

\subsection{Previous work}

Let $\mP$ be a set of points and $\mL$ a set of lines in $\R^2$ with $|\mP|=m,|\mL|=n$. The well-known Szemerédi-Trotter theorem \cite[Theorem 1]{Szemeredi1983extremal} states that
\begin{equation}
\I(\mP,\mL)=O(m^{2/3}n^{2/3}+m+n),
\end{equation}
where 
$$\I(\mP,\mL)=|\{(p,l): p\in\mP,l\in\mL,p\text{ lies on }l|$$ 
denotes the number of incidences. This bound is sharp and cannot be improved. 

The point-line incidences have also been studied in other settings. The Szemerédi-Trotter theorem in the complex plane $\C^2$ is proved by Tóth \cite[Theorem 1]{toth2015Szemeredi}, which turns out to have the same expression as the real case above. Over finite fields, Bourgain, Katz, and Tao \cite[Theorem 6.2]{bourgain2004sum} proved an incidence bound of Szemerédi–Trotter type for $\F_q$ using sum–product estimates, and Vinh \cite[Theorem 3]{vinh2011Szemeredi} later obtained another bound via spectral graph theory, which becomes stronger when the number of points and lines is large. See also \cite{rudnev2018number,stevens2017improved,grosu2014f,iosevich2023improved} for further developments in this direction. Incidence geometry has deep connections with other areas, including geometric measure theory, additive combinatorics, and harmonic analysis. One may turn to Dvir \cite{dvir2012incidence} for a summary.

\subsection{Main results}

Motivated by these developments, it is natural to study point-line incidences over more general fields, such as the $p$-adic field $\Q_p$ and its extensions, or other transcendental extensions over $\Q$. The result of this paper covers all these cases. In fact, we establish the Szemerédi–Trotter theorem over arbitrary fields of characteristic zero. From now on, $K$ is always a field with $\char(K)=0$. Our main result is the following.

\begin{thm}\label{thm: incidence bound_main theorem}
Let $\mP$ be a set of points and $\mL$ a set of lines in $\K^2$ with $|\mP|=m,|\mL|=n$. Then we have
\begin{equation}\label{eq: TS incidence bound}
\I(\mP,\mL)=O(m^{2/3}n^{2/3}+m+n),
\end{equation}
and this bound is sharp.
\end{thm}

The result above strengthens a theorem of Stevens and De Zeeuw \cite[Theorem 3]{stevens2017improved}. The key step to prove \Cref{thm: incidence bound_main theorem} is to apply the following Baby Lefschetz principle, which first appeared in the appendix of \cite{lefschetz2005algebraic}. See also Tao's blog post \cite[Proposition 4]{tao_lefschetz} for related discussions.

\begin{lemma}\label{lem: Baby Lefschetz principle}
(Baby Lefschetz principle) Let $F$ be a field of characteristic zero that is finitely generated over $\Q$. Then there exists a field isomorphism $\phi:F\rightarrow\phi(F)$ from $F$ to a subfield $\phi(F)$ of $\C$.
\end{lemma}

Informally, \Cref{lem: Baby Lefschetz principle} allows us to embed any finite configuration of points and lines in $K^2$ into $\C^2$ without changing the incidence relations; see \Cref{prop: complex plane embedding}. Since the complex plane case has already been proved in \cite[Theorem 1]{toth2015Szemeredi}, \Cref{thm: incidence bound_main theorem} follows.

The Baby Lefschetz principle also has other applications, matching the incidence behavior in the field $K$ to that in the complex case. We write $A\lesssim B$ if $A\le CB$ for some constant $C$. A proposition that might be useful is the following.

\begin{prop}\label{prop: estimate of L_n}
Let $\mP$ and $\mL$ be finite sets of points and lines in $K^2$. For all $n\ge 2$, denote by $\mL_n$ the set of lines in $\mL$ containing at least $n$ points in $\mP$. Then
$$|\mL_n|\lesssim\frac{|\mP|^2}{n^3}+\frac{|\mP|}{n}.$$
\end{prop}

The following theorem is the analog of Beck's theorem (see \cite[Theorem 3.1]{beck1983lattice}) over $K$.

\begin{thm}\label{thm: Beck's theorem for K}
Let $\mP$ be a finite set of points in $K^2$, and let $\mL$ be the set of lines that contain at least two points of $\mP$. Then, at least one of the following is true:
\begin{enumerate}
\item There exists a line in $\mL$ that contains $\gtrsim|\mP|$ points of $\mP$.
\item $|\mL|\gtrsim|\mP|^2$.
\end{enumerate}
\end{thm}

We can also prove incidence results between points and other kinds of algebraic varieties in $K$. One may move to \Cref{sec: more general incidences} for further discussions around this.

For a finite set $A\subset K$, the set of pairwise sums and products formed by elements of $A$ are given by
$$A+A=\{a+b\mid a,b\in A\},A\cdot A=\{ab\mid a,b\in A\}$$
respectively. Applying the Baby Leschetz principle to the sum-product estimate over the complex field in \cite[Theorem 1.1]{konyagin2013new}, we obtain the following theorem:

\begin{thm}\label{thm: sum-product estimate}
(Erd\H{o}s–Szemerédi theorem) For a finite subset $A\subset K$, we have
$$\max\{|A+A|,|A\cdot A|\}\gtrsim\frac{|A|^{4/3}}{(\log |A|)^{1/3}}.$$
\end{thm}

There have been several attempts to improve sum–product estimates over certain specific fields of characteristic zero, beyond the classical real or complex settings. For example, Bloom and Jones \cite[Theorem 1.3]{bloom2014sum} obtained a sum–product estimate over $\Q_p$, and Croot and Hart \cite[Theorem 1]{croot2010sums} proved an analogue over $\C(x)$\footnote{Strictly speaking, Croot and Hart worked over the ring $\C[x]$. Since any finite subset of $\C(x)$ can be multiplied by a common denominator to obtain a subset of $\C[x]$ without affecting the sum-product properties, their result extends to $\C(x)$ in the sense used here.}. While these results rely on the specific algebraic or analytic structures of the underlying field, our observation improves these lower bounds by reducing the problem to the complex setting. This suggests that, for fields of characteristic zero, a viable approach of related problems is to bypass the special structure of the field altogether and instead relate the problem directly to the complex numbers. The justification for this reduction lies in a general principle from model theory (see \cite{barwise1969lefschetz}), which asserts that all first-order sentences over algebraically closed fields of characteristic zero are equivalent. The baby Lefschetz principle is in fact a special case of this phenomenon. Since both incidence statements involving finitely many objects and sum–product estimates for finite sets can be formulated as first-order statements, they may always be transferred to the complex field $\C$, the algebraically closed field most familiar in analysis and algebraic geometry. We hope this connection will lead to further applications: once an incidence or sum–product result is established over $\C$, an identical statement automatically follows over any field $K$ of characteristic zero.

\subsection{Plan of the paper}

In \Cref{sec: Proof of the main results}, we apply \Cref{lem: Baby Lefschetz principle} to prove the results in \Cref{sec: Introduction}. In \Cref{sec: more general incidences}, we extend our method to incidence problems involving more general algebraic varieties and to higher-dimensional analogues.

\section{Proof of the main results}\label{sec: Proof of the main results}

In this section, we prove the results stated in \Cref{sec: Introduction}. \Cref{lem: Baby Lefschetz principle} leads to the following proposition, which intuitively states that we can embed points and lines into the complex plane while preserving the incidence relations.

\begin{prop}\label{prop: complex plane embedding}
Let $\mP$ be a set of points and $\mL$ a set of lines in $\K^2$ with $|\mP|=m,|\mL|=n$. Then there exists an injective map $\phi_\mP$ from $\mP$ to points in $\C^2$, and an injective map $\phi_{\mL}$ from $\mL$ to lines in $\C^2$, such that for all $p\in \mP$ and $l\in\mL$, $p$ lies on $l$ if and only if $\phi_{\mP}(p)$ lies on $\phi_{\mL}(l)$.
\end{prop}

\begin{proof}
We denote
$$\mP=\{(x_1,y_1),\ldots,(x_m,y_m)\},\mL=\{a_1x+b_1y+c_1=0,\ldots,a_nx+b_ny+c_n=0\},$$
where $x_i,y_i\in K$ for all $1\le i\le m$, and $a_i,b_i,c_i\in K$ for all $1\le i\le n$. Let 
$$F=\Q(x_1,y_1,\ldots,x_m,y_m,a_1,b_1,c_1,\ldots,a_n,b_n,c_n).$$
By \Cref{lem: Baby Lefschetz principle}, there exists a field isomorphism $\phi:F\rightarrow\phi(F)$ from $F$ to a subfield $\phi(F)$ of $\C$. Let
$$\phi_\mP:(x_i,y_i)\mapsto(\phi(x_i),\phi(y_i)),\quad 1\le i\le m,$$
and 
$$\phi_\mL:a_ix+b_iy+c_i=0\mapsto \phi(a_i)x+\phi(b_i)y+\phi(c_i)=0,\quad 1\le i\le n.$$
Then for all $1\le i\le m$ and $1\le j\le n$, we have
\begin{align}
\begin{split}
(x_i,y_i)\text{ lies on }a_jx+b_jy+c_j=0&\Leftrightarrow a_jx_i+b_jy_i+c_j=0\\
&\Leftrightarrow\phi(a_jx_i+b_jy_i+c_j)=0\\
&\Leftrightarrow(\phi(x_i),\phi(y_i))\text{ lies on }\phi(a_j)x+\phi(b_j)y+\phi(c_j)=0\\
&\Leftrightarrow\phi_\mP(x_i,y_i)\text{ lies on }\phi_{\mL}(a_jx+b_jy+c_j)=0,
\end{split}
\end{align}
which ends the proof.
\end{proof}

\begin{rmk}
Note that the field $K$ might be very large that we cannot embed the entire plane $K^2$ into $\C^2$. However, for our purposes, it suffices to embed only the finitely many points and lines under consideration into $\C^2$.
\end{rmk}

Now we can start our proof of \Cref{thm: incidence bound_main theorem}.
\begin{proof}[Proof of \Cref{thm: incidence bound_main theorem}]
On the one hand, take $\phi_{\mP},\phi_{\mL}$ as in \Cref{prop: complex plane embedding}. Then we have
$$\I(\mP,\mL)=\I(\phi_{\mP}(\mP),\phi_{\mL}(\mL))=O(m^{2/3}n^{2/3}+m+n),$$
where the second equality follows from the Szemerédi-Trotter theorem in the complex plane, as proved in \cite[Theorem 1]{toth2015Szemeredi}.

On the other hand, recall the sharp construction from Elekes \cite[Example 1.16]{elekes2001sums} for the real plane. In this construction, all points have integer coordinates and all lines have integer coefficients. Therefore, the same example gives a sharp construction over $K$ as well.

\end{proof}

The strategy of the proof of \Cref{prop: estimate of L_n} and \Cref{thm: Beck's theorem for K} is similar to the proof of \Cref{thm: incidence bound_main theorem}. We embed the points and lines into the complex plane $\C^2$, and then the incidence in $K^2$ follows directly from the results in $\C^2$ that are already proven in \cite{toth2015Szemeredi}.

\begin{proof}[Proof of \Cref{prop: estimate of L_n}]
Take $\phi_{\mP},\phi_{\mL}$ as in \Cref{prop: complex plane embedding}. For all $n\ge 2$, let $\phi_{\mL}(\mL_n)$ denote the set of lines in $\phi_{\mL}(\mL)$ that contain at least $n$ points of $\phi_{\mP}(\mP)$. By \cite[Theorem 2]{toth2015Szemeredi}, we have
\begin{align}
\begin{split}
|\mL_n|&=|\phi_{\mL}(\mL_n)|\\
&\lesssim\frac{|\phi_{\mP}(\mP)|^2}{n^3}+\frac{|\phi_{\mP}(\mP)|}{n}\\
&=\frac{|\mP|^2}{n^3}+\frac{|\mP|}{n},
\end{split}
\end{align}
which completes the proof.
\end{proof}

\begin{proof}[Proof of \Cref{thm: Beck's theorem for K}]
Take $\phi_{\mP},\phi_{\mL}$ as in \Cref{prop: complex plane embedding}. Applying \cite[Corollary 3]{toth2015Szemeredi}, at least one of the following holds:
\begin{enumerate}
\item There exists a line in $\phi_{\mL}(\mL)$ that contains $\gtrsim|\mP|$ points of $\phi_{\mP}(\mP)$. In this case, there exists a line in $\mL$ that contains $\gtrsim|\mP|$ points of $\mP$.
\item $|\phi_{\mL}(\mL)|\gtrsim|\phi_{\mP}(\mP)|^2$, which implies $|\mL|\gtrsim|\mP|^2$.
\end{enumerate}
\end{proof}

To prove the sum-product estimate in \Cref{thm: sum-product estimate}, we aim to find a map $\phi$ that embeds $A\subset K$ into $\C$.

\begin{proof}[Proof of \Cref{thm: sum-product estimate}]
Denote $A=\{a_1,\ldots,a_n\}$, where $n=|A|$. Let $F=\Q(a_1,\ldots,a_n)$. By \Cref{lem: Baby Lefschetz principle}, there exists a field isomorphism $\phi:F\rightarrow\phi(F)$ from $F$ to a subfield $\phi(F)$ of $\C$. Then we have
\begin{align}
\begin{split}
\max\{|A+A|,|A\cdot A|\}&=\max\{|\phi(A)+\phi(A)|,|\phi(A)\cdot\phi(A)|\}\\
&\ge C\frac{|\phi(A)|^{4/3}}{(\log |\phi(A)|)^{1/3}}\\
&=C\frac{|A|^{4/3}}{(\log |A|)^{1/3}},
\end{split}
\end{align}
where $C>0$ is the same implicit constant in \cite[Theorem 1.1]{konyagin2013new}.
\end{proof}

\begin{rmk}
Currently, the sum-product estimate in \Cref{thm: sum-product estimate} gives only the bound  $\gtrsim\frac{|A|^{4/3}}{(\log |A|)^{1/3}}$, as this is the best known result over the complex numbers at the time of writing. Nevertheless, by applying exactly the same method, any future improvement of the complex sum–product bound would automatically lead to a corresponding improvement of our lower bound over $K$. In particular, Erd\H{o}s \cite{erdos1976some} conjectured that 
$$\max\{|A+A|,|A\cdot A|\}\gtrsim |A|^{2-o(1)}.$$
It is therefore natural to make the same conjecture over $K$. Moreover, this bound would be sharp, since the construction in \cite{erdos1983sums} involves only integer sets and thus also applies to subsets of $K$.
\end{rmk}

\section{More general incidences}\label{sec: more general incidences}

In addition to the point-line incidences studied in \Cref{sec: Proof of the main results}, one may also consider the incidence between points and other geometric objects, including circles, planes, and other types of algebraic varieties. Let $\mP$ be a set of points, and $\mathcal{V}$ be a set of geometric objects. Denote by
$$\I(\mP,\mathcal{V})=|\{(p,V):p\in\mP,V\in\mathcal{V},p\text{ lies on }V\}|$$
the incidences between $\mP$ and $\mathcal{V}$. In this section, we discuss general incidence results in $K$, which are all first-order statements and therefore can be deduced by applying \Cref{lem: Baby Lefschetz principle} to known results in the complex field.

Let $\mP$ be a set of points and $\mathcal{C}$ be a set of curves in $K^2$. Following the definition on the first page of \cite{sheffer2018point}, we say that $(\mP,\mathcal{C})$ \emph{has $k$ degrees of freedom and multiplicity type $s$} if 
\begin{enumerate}
\item For any $\mP'\subset\mP$ of size $k$, there are at most $s$ curves in $\mathcal{C}$ that contain $\mP'$.
\item Any two curves in $\mathcal{C}$ intersect at most $s$ points of $\mP$.
\end{enumerate}

Our result is as follows, which is the analog of \cite[Theorem 1.3]{sheffer2018point} to $K$.

\begin{thm}
Let $k\ge 1,D\ge 1,s\ge 1$, and $\epsilon>0$. Let $\mP\subset K^2$ be a set of $m$ points, and $\mathcal{C}$ be a set of $n$ algebraic curves over $K$, each of degree at most $D$. Suppose that $(\mP,\mathcal{C})$ has $k$ degrees of freedom and multiplicity type $s$. Then
$$\mathcal{I}(\mP,\mathcal{C})\le C(m^{\frac{k}{2k-1}+\epsilon}n^{\frac{2k-2}{2k-1}}+m+n),$$
where the definition of $C=C(\epsilon,D,s,k)$ is the same as in \cite[Theorem 1.3]{sheffer2018point}.
\end{thm}

\begin{proof}
Let $F$ be the smallest field over $\Q$ that contains all the coordinates of the points in $\mP$ and the coefficients of the curves in $\mathcal{C}$. Then $F$ is finitely generated. By \Cref{lem: Baby Lefschetz principle}, there exists a field isomorphism $\phi:F\rightarrow\phi(F)$ from $F$ to a subfield $\phi(F)$ of $\C$. Define injective maps $\phi_{\mP}$ and $\phi_{\mC}$, analogous to those in \Cref{prop: complex plane embedding}, that embed $(\mP,\mathcal{C})$ into the complex plane. Then $(\phi_{\mP}(\mP),\phi_{\mathcal{C}}(\mathcal{C}))$ also has $k$ degrees of freedom and multiplicity type $s$. Applying \cite[Theorem 1.3]{sheffer2018point}, we have
$$\I(\mP,\mathcal{C})=\I(\phi_{\mP}(\mP),\phi_{\mathcal{C}}(\mathcal{C}))\le C(m^{\frac{k}{2k-1}+\epsilon}n^{\frac{2k-2}{2k-1}}+m+n).$$
\end{proof}

We can also study incidences in higher-dimensional spaces. Let $d\ge 1$ be an integer and $\mP\subset K^d$ be a set of $m$ points. Let $n\ge 2$, and $\mL_n(\mP)$ be the set of lines that are incident to at least $n$ points from $\mP$. The following theorem originates from \cite[Theorem 1.3]{zahl2016note}, which intuitively states that if a collection of points in $K^d$ gives many $n$-rich lines, then a positive proportion of these points must lie on a common $(d-1)$-flat.

\begin{thm}
Let $d\ge 1$ and $\epsilon>0$. Let $\mP\subset K^d$ be a set of $m$ points, and let $n\ge 2$. Suppose that
$$|\mL_n(\mP)|>C_{d,\epsilon}\cdot\alpha\cdot\frac{n^{2+\epsilon}}{r^{d+1}}$$
for some $\alpha\ge 1$. Then there exists a subset $\mP'\subset \mP$ with $|\mP'|\ge c_{d,\epsilon}\cdot\alpha\cdot\frac{n^{2+\epsilon}}{r^{d+1}}$, which is contained in a $(d-1)$-flat. Here the constants $c_{d,\epsilon},C_{d,\epsilon}$ are the same as in \cite[Theorem 1.3]{zahl2016note}.
\end{thm}

\begin{proof}
Let $F$ be the smallest field over $\Q$ containing all the coordinates of the points in $\mP$ and the coefficients of the lines in $\mL_n$. Then $F$ is finitely generated. By \Cref{lem: Baby Lefschetz principle}, there exists a field isomorphism $\phi:F\rightarrow\phi(F)$ from $F$ to a subfield $\phi(F)$ of $\C$. Define injective maps $\phi_{\mP},\phi_{\mL_n}$ in a way similar to \Cref{prop: complex plane embedding}, which embed $(\mP,\mL_n)$ into $\C^d$. By \cite[Theorem 1.3]{zahl2016note}, the result follows.
\end{proof}

\begin{rmk}
In fact, the cheap Dvir-Gopi version (see \cite[Corollary 1.1]{zahl2016note}) also holds for any field $K$ of characteristic zero. The proof is essentially the same as above, so we omit it.
\end{rmk}

\bibliographystyle{plain}
\bibliography{references.bib}

\end{document}